\documentclass[12pt]{amsart}
\usepackage{pdfsync}
\usepackage{amsmath,amsthm,amsfonts,amssymb}
\usepackage{amsmath, amsthm, amssymb, amscd, amsfonts}
\usepackage[all]{xy}
\usepackage{amssymb, amsthm, amsmath}
\usepackage[english]{babel}
\usepackage{amsfonts}
\numberwithin{equation}{subsection}

\DeclareMathOperator{\id}{id}

\newcommand{\ra}{\rightarrow}

\newcommand{\ot}{\otimes}
\newcommand{\mtc}{\mathcal}

\newcommand{\Lam}{\Lambda}

\newcommand{\al}{\alpha}
\newcommand{\eps}{\epsilon}

\newcommand{\D}{\Delta}

\newtheorem{lemma}[equation]{Lemma}
\newtheorem{thm}[equation]{Theorem}
\newtheorem{prop}[equation]{Proposition}
\newtheorem{defn}[equation]{Definition}
\newtheorem{cor}[equation]{Corollary}
\newtheorem{rem}[equation]{Remark}

\newcommand{\dw}{\downarrow}
\newcommand{\uw}{\uparrow}

\newcommand{\ch}{\chi}
\newcommand{\mtr}{\mathrm}
\newcommand{\bn}{\begin}

\newcommand{\ncm}{\newcommand}
\ncm{\np}{\newpage}
\ncm{\ebl}{\end{thebibliography}}
\ncm{\bbl}{\begin{thebibliography}}
\ncm{\chd}{_{ _{\ch}}}
\ncm{\ald}{_{ _{\al}}}
\ncm{\cP}{\mathcal{P}}
\ncm{\ei}{e_i}
\ncm{\eij}{e_{i,\;j}}
\ncm{\bt}{\begin{thm}}
\ncm{\bdef}{\begin{defn}}
\ncm{\edf}{\end{defn}}
\ncm{\et}{\end{thm}}
\ncm{\bc}{\begin{cor}}
\ncm{\bl}{\begin{lemma}}
\ncm{\el}{\end{lemma}}
\ncm{\bpf}{\begin{proof}}
\ncm{\epf}{\end{proof}}
\ncm{\ec}{\end{cor}}
\ncm{\er}{\end{rem}}
\ncm{\br}{\begin{rem}}
\ncm{\bp}{\begin{prop}}
\ncm{\ep}{\end{prop}}
\ncm{\bd}{\begin{document}}
\ncm{\ed}{\end{document}}
\ncm{\beq}{\begin{equation}}
\ncm{\beqn}{\begin{equation*}}
\ncm{\eeq}{\end{equation}}
\ncm{\eeqn}{\end{equation*}}
\ncm{\bea}{\begin{eqnarray}}
\ncm{\eea}{\end{eqnarray}}
\ncm{\beanon}{\begin{eqnarray*}}
\ncm{\eeanon}{\end{eqnarray*}}
\title[Semisimple Hopf algebras]
{Kernels of representations and coideal subalgebras of Hopf algebras}
\author{Sebastian  Burciu}
\address{Inst.\ of Math.\ ``Simion Stoilow" of the Romanian Academy
P.O. Box 1-764, RO-014700, Bucharest, Romania}
\email{sebastian.burciu@imar.ro}
\date{\today}
\begin{document}
\maketitle
\bn{abstract}
We define left and right kernels of representations of Hopf algebras. In the case of group algebras, left and right kernels coincide and they are the usual kernels of modules. In the general case we show that these kernels coincide with the categorical left and right Hopf kernels of morphisms of Hopf algebras defined in \cite{AD}. Brauer's theorem for kernels over group algebras is generalized to Hopf algebras.
\end{abstract}
\section{Introduction}
If $G$ is a finite group and $M$ is a representation of $G$ which affords the character $\ch$ then the kernel of the representation $M$ is defined as the set of all $g \in G$ which act as the identity on $M$.  It is easy to see that $\mtr{ker}_G(\ch)$ is a normal subgroup of $G$. The notion of kernel of a representation of a group is a very important notion in studying the representations of groups (see \cite{Is}). For example, a classical result of Brauer in group theory states that over an algebraically closed field $k$ of characteristic zero if $\ch$ is a faithful
character of $G$ then any other irreducible character of $G$ is a direct summand
in some tensor power of $\ch$.

In \cite{Bker} the author introduced a similar notion for the kernel
of representations of any semisimple Hopf algebra. The notion uses the
fact that the exponent of a semisimple Hopf algebra is finite \cite{EG'}.
In this paper we will extend the notion of kernel to arbitrary Hopf algebras,
not necessarily semisimple. The aforementioned Brauer's theorem for groups was extended in \cite{Bker} for semisimple Hopf algebras $A$ and for those representations of $A$ whose characters are central in $A^*$. With the help of this new notion of kernel, we extend this theorem to an arbitrary representation whose character is not necessarily central in the dual Hopf algebra, see Theorem \ref{Brauer}. We also prove a version of Brauer's theorem for arbitrary Hopf algebras, giving a new insight to the main results from \cite{PQ} and \cite{Riburn}.

Note that the group algebra $kG$ is cocommutative. Lack of cocommutativity suggests the introduction of left and right kernels of modules. It will also be shown that these left and right kernels coincide with the left and respectively right categorical kernels of morphisms of Hopf algebras that has been introduced in \cite{AD} and extensively studied in \cite{Andr}. This coincidence proven in Theorem \ref{catchr} suggests also how to generalize the notion of kernel of modules to the non semisimple case.

Recently it was proven in \cite{BJ} that Hopf subalgebras are normal if and only if they are depth two subalgebras. We extend this result to coideal subalgebras. Moreover we show that in this situation, depth two and normality in the sense defined by Rieffel in \cite{Ri} also coincide. Normal subalgebras as defined by Rieffel were recently revised in Section 4 of \cite{BKK}.

The paper is organized as follows.
The second section introduces the notion of left and right kernels of modules and proves their coincidence with the categorical kernels introduced in \cite{AD}. A general version of Brauer's theorem is also stated in this section. The third section is concerned with
depth two coideal subalgebras of a Hopf algebra. We prove a left coideal subalgebra is
right depth two subalgebra if and only if it is a left normal coideal subalgebra. Section \ref{ss} considers coideal subalgebras of semisimple Hopf algebra. An extension of Brauer's theorem in this case
is also considered.

The Hopf algebra notations from \cite{Montg} are used in this paper but we drop the sigma symbol from Sweedler's notation of comultiplication. Recall a left coideal subalgebra $S$ of $A$ is a subalgebra of $A$ with $\D(S)\subset A\ot S$. We say that a left coideal subalgebra $S$ of $A$ is left normal if and only if it is invariant
under the left the adjoint action of $A$, i.e $a_1xS(a_2) \in S$ for all $a \in A$ and $x \in S$. Such a coideal subalgebra will be called a left normal coideal subalgebra. Also it follows from \cite{Sk} that if $A$ is finite dimensional then $A$ is free as left or right $S$-module.

For a subalgebra $B$ of a Hopf algebra $A$ we denote by $B^+$ the ideal of $B$ defined
by $\omega(B)=B^+=A^+\cap B$. As usual, $A^+=\mtr{ker}\;\eps$ is the kernel of the
count $\eps$ of $A$.
\section{Left and right kernels}\label{lrk}
Let $A$ be an arbitrary Hopf algebra over a field $k$ and $M$ be an $A$-module.
We say that an element $a \in A$ acts trivially on $M$ if and only
$am=\eps(a)m$ for all $a \in A$. If $S$ is a subalgebra of $A$ we say that $S$ acts trivially
on $M$ if each element of $S$ acts trivially on $M$.

If $M$ is an $A$-module then define

\begin{equation}\label{A}
    {\bf{{\bf \mtc{A}}}_{ _{M}}}=\{a\in A\;|\; am=\eps(a)m\;\text{for all} \;m\in M\}
\end{equation}

It is easy to verify that ${\bf \mtc{S}}_{ _{M}}$ is a subalgebra of $A$, thus the largest subalgebra of $A$ which acts trivially on $M$. If $A$ is semisimple and $M$ a simple $A$-module then $\dim_k({\bf \mtc{S}}_{ _{M}})=\dim_k(A)-\dim_k(M)^2$.
\subsection{Definition of left kernels} Let $M$ be an $A$-module and $\mtr{L}_{ _M}:A \ot M\ra A\ot M$ be the linear operator given by $a\ot m\mapsto a_1 \ot a_2m $. Let $\mtr{LKer}_{ _M}$ be the largest subspace $B$ of $A$ such that $\mtr{L}_{ _M}|_{ _{B\ot M}}=\mtr{id}$. Then  $\mtr{LKer}_{ _M}$ is called the left kernel of $M$. Thus

\begin{equation}\label{L}
\mtr{LKer}_{ _{M}}=\{a \in A|\; a_1\ot a_2m=a\ot m,\;\;\text{for all}\;\; m\in M\}
\end{equation}

If $A$ is a semisimple Hopf algebra and $M$ a module with character $\ch$ we also write $\mtr{LKer}_{ _{\ch}}$ instead of $\mtr{LKer}_{ _M}$.
\bn{rem}
If $A=kG$ and $M$ a $G$-module then $\mtr{LKer}_{ _M}$ is the group algebra $k[\mtr{Ker}_G(M)]$ of the usual kernel of $M$.
\end{rem}

\begin{prop}\label{charatcerization}
Let $A$ be a Hopf algebra and $M$ be an $A$-module. Then $\mtr{LKer}_{ _{M}}$ is a left normal coideal subalgebra of $A$.
\end{prop}

\bn{proof}
Applying $\eps\ot \mtr{Id}$ to the defining relation \ref{L} of $\mtr{LKer}_{ _{M}}$ it follows that $\mtr{LKer}_{ _{M}}\subset {\bf \mtc{S}}_{ _{M}}$. It is easy to verify that $\mtr{LKer}_{ _{M}}$ is an algebra.

Next we will show that $\D(\mtr{LKer}_{ _{M}})\subset A  \ot \mtr{LKer}_{ _{M}}$.
Let $a\in \mtr{LKer}_{ _{M}}$ and $\D(a)=\sum_{i=1}^nh_i \ot x_i$ where $h_i$ is a basis of $A$. It will be shown that $x_i \in \mtr{LKer}_{ _{M}}$ for all $1 \leq i \leq n$. Indeed $\D^2(a)=\sum_{i=1}^nh_i\ot \D(x_i)$. Since $a_1\ot a_2m=a\ot m$ it follows that $a_1\ot a_2\ot a_3m=a_1\ot a_2\ot m$ for all $m \in M$.
Thus $\sum_{i=1}^nh_i \ot (x_i)_1\ot (x_i)_2m=\sum_{i=1}^nh_i \ot x_i\ot m$ which implies that $(x_i)_1\ot (x_i)_2m=x_i\ot m$ for all $m \in M$ and $1 \leq i \leq n$.

Now we will show invariance under the left adjoint action of $A$. Let $h \in A$ and $a\in \mtr{LKer}_{ _{M}}$. Then $h_1aS(h_2)\in \mtr{LKer}_{ _{M}}$ since
\bn{eqnarray*}
(h_1aS(h_2))_1\ot (h_1aS(h_2))_2m & \! = \! & h_1a_1S(h_4)\ot h_2a_2S(h_3)m \\ & = & h_1aS(h_4)\ot h_2S(h_3)m \\ & = & h_1aS(h_2)\ot m
\end{eqnarray*}
\end{proof}

\bp\label{seseq}
Let $A$ be a Hopf algebra and
$$\begin{CD}0 @> >> N @> i>>  M @> p >>P  @> >> 0 \end{CD}$$

be a short exact sequence of modules. Then
$$\mtr{LKer}_ {_{M}} \subseteq \mtr{LKer}_ {_{
N}}\cap \mtr{LKer}_ {_{P}}.$$

\ep

\bpf
Clearly $\mtr{LKer}_ {_{M}} \subset \mtr{LKer}_ {_{
N}}$. Next we will show  $\mtr{LKer}_{ _M}\subset \mtr{LKer}_{ _P}$.
Indeed if $a \in \mtr{LKer}_{ _M}$ then $a_1\ot a_2m=a\ot m$ for all $m \in M$ and applying $\mtr{id} \ot p$ it follows that $a \in \mtr{LKer}_{ _P}$ since $p$ is surjective.
\epf
Recall that the exponent of $A$ is the smallest positive number $m >0$ such that $a^{[m]}=\eps(a)1$ for all
$a \in A$. The generalized power $a^{[m]}$ is defined by
$a^{[m]}=a_1a_2...a_m$.

\bn{rem}\end{rem}Let $A$ be a Hopf algebra  and $M$ be an $A$-module. Then it is easy to check that $L_{M}^s(a\ot m)=a_1\ot a_2^{[s]}m$ for all $s \geq 1$. Thus if $A$ is of a finite exponent $m$ then $L_M^m=\id$ on $A\ot M$.

\bn{thm}
Let $A$ be a Hopf algebra of finite exponent over an algebraically closed field $k$ and let
$$\begin{CD}0 @> >> N @> i>>  M @> p >>P  @> >> 0 \end{CD}$$

be a short exact sequence of finite dimensional modules. Then
$$\mtr{LKer}_ {_{M}} = \mtr{LKer}_ {_{
N}}\cap \mtr{LKer}_ {_{P}}.$$
\end{thm}
\bpf
Let $L:=\mtr{LKer}_ {_{ N}}\cap \mtr{LKer}_ {_{P}}$ and consider the exact sequence

$$\begin{CD}0 @> >>L\ot N @> 1\ot i>> L\ot  M @> 1\ot p >>L\ot P  @> >> 0 \end{CD}.$$

If $l \in L$ and $m \in M$ then $$(1\ot p)(L_{ _M}(l\ot m)-l \ot m)= l_1\ot p(l_2m)-l\ot p(m)=l_1\ot l_2p(m)-l\ot p(m)=0$$
since $l\in \mtr{LKer}_{ _P}$. Thus $L_{ _M}(l\ot m)-l \ot m\in L\ot N$. It follows that $L_{ _N}(L_{ _M}(l\ot m)-l \ot m))=L_{ _M}(l\ot m)-l \ot m$. This also implies that $L_M^2(l \ot m)-L_M(l \ot m)=L_{ _M}(l \ot m)-l \ot m$. Thus $(L_{ _M}-\mtr{id})^2=0$ on $L\ot M$. On the other hand since $A$ has finite exponent one has $L_{ _M}^m=\mtr{id}$ on $L\ot M$ by the previous remark. Since $k$ is algebraically closed it follows that $L_{ _M}=\id$ on $L\ot M$, i.e $L \subset \mtr{LKer}_ {_{M}}$.
\epf

\bn{rem}\end{rem}
Among the Hopf algebras with finite exponent we recall semisimple algebras and group algebras of finite groups over a field of characteristic diving the order of the group. The fact that the exponent of a finite dimensional semisimple Hopf algebra is finite was proven in \cite{EG'}. In the same paper it is also shown that in this case the exponent divides the
third power of the dimension of $A$.

\bn{rem}\label{bij}\end{rem}
Let $S_{ _M}$ the inverse operator of $L_{ _M}$. Then $S_{ _M}$ is given by $S_{ _M}(a\ot m)=a_1\ot S(a_2)m$ for all $a \in A$ and $m \in M$. It follows that $S_{ _M}|_{ _{\mtr{LKer}_{ _{M}}\ot A}}=\mtr{id}$ and thus one also has that: $$\mtr{LKer}_{ _{M}}=\{a \in A|\; a_1\ot S(a_2)m=a\ot m,\;\;\text{for all}\;\;m\in M\}$$

\bp\label{largestcoidsubalg}
Let $A$ be a Hopf algebra and $M$ be an $A$-module. Then $\mtr{LKer}_{ _M}$ is the largest coideal subalgebra of $A$ that acts trivially on $M$.
\ep

\bpf
Proposition \ref{charatcerization} implies that $\mtr{LKer}_{ _M}$ acts trivially on $M$. On the other hand if $\D(S)\subset A\ot S$ and $S$ acts trivially on $M$ then $s_1\ot  s_2m=s_1\ot\eps(s_2)m=s\ot m$, which shows that $S\subset \mtr{LKer}_{ _M}$.
\epf

\bn{lemma}\label{dual}
Let $A$ be a Hopf algebra with bijective antipode and $M$ be a finite dimensional module over $A$. Then $\mtr{LKer}_{ _{M}}= \mtr{LKer}_{ _{M^*}}$.
\end{lemma}

\bn{proof}
For the operator $L_{ _{M^*}}:A \ot M^*\ra A\ot M^*$ one has that $\mtr{LKer}_{ _{M^*}}=\{a \in A|\; a_1\ot a_2f=a\ot f\}$ for all $f \in M^*$. This is equivalent to $a_1f(S(a_2)m)=af(m)$ for all $m \in M$ and $f\in M^*$. Therefore $a \in \mtr{LKer}_{ _{M^*}}$ if and only if $a_1\ot S(a_2)m=a\ot m$ for all $m \in M$ which implies by Remark \ref{bij} that $\mtr{LKer}_{ _{M^*}}=\mtr{LKer}_{ _{M}}$.
\end{proof}
\subsection{Definition of right kernel}
Let $M$ be an $A$-module and \\$\mtr{R}_{ _M}:A \ot M\ra A\ot M$ given by $a\ot m\mapsto a_2 \ot a_1m $. Let $\mtr{RKer}_{ _{M}}$ be the largest subspace $B$ of $A$ such that $\mtr{R}_{ _M}|_{ _{B\ot M}}=\mtr{id}$. Thus

\begin{equation}\label{R}
\mtr{RKer}_{ _{M}}=\{a \in A|\; a_2\ot a_1m=a\ot m,\;\;\text{for all}\;\;m\in M\}
\end{equation}
Suppose that $A$ has a bijective antipode $S$ with inverse $S^{-1}$.
Then the inverse operator $U_{ _M}$ of $\mtr{R}_{ _M}$ is given by $U_{ _M}(a\ot m)=a_2\ot S^{-1}(a_1)m$ for all $a \in A$ and $m \in M$. It follows that $U_{ _M}|_{ _{\mtr{RKer}_{ _{M}}\ot A}}=\mtr{id}$ and thus one has also that

\begin{equation}\label{R'}
    \mtr{RKer}_{ _{M}}=\{a \in A|\; a_2\ot S^{-1}(a_1)m=a\ot m,\;\;\text{for all}\;\;m\in M\}
\end{equation}
 It is also easy to check that $\D(\mtr{RKer}_{ _{M}})\subset \mtr{RKer}_{ _{M}} \ot  A $.
\bn{rem}\end{rem}
1)Suppose that the Hopf algebra $A$ has a bijective antipode and let $M$ be an $A$-module. Then $\mtr{RKer}_{ _{M}}=S(\mtr{LKer}_{ _{M}})$.
Indeed if $a \in \mtr{LKer}_{ _{M}}$ then
\bn{eqnarray*}
S(a)_2\ot S(a)_1m & \! = \! & S(a_1)\ot S(a_2)m  =   S(a)\ot m
\end{eqnarray*}
which shows that $S(\mtr{LKer}_{ _{M}}) \subset \mtr{RKer}_{ _{M}}$. Similarly it can be checked that if $a \in \mtr{RKer}_{ _{M}}$ then $S^{-1}a \in \mtr{LKer}_{ _{M}}$. Thus $\mtr{RKer}_{ _{M}}=S(\mtr{LKer}_{ _{M}})$.
\\
2)Applying $\eps \ot\mtr{Id}$ to the relation \ref{R} it also follows that $\mtr{RKer}_{ _{M}}\subset {\bf \mtc{S}}_{ _{M}}$.

\subsection{Description as categorical kernels}
\bn{lemma}\label{power}
If $M$ and $N$ are two $A$-modules then $\mtr{LKer}_{ _{M}}\cap \mtr{LKer}_{ _{N}} \subset \mtr{LKer}_{ _{M\otimes N}}$. In particular $\mtr{LKer}_{ _{M}} \subset \mtr{LKer}_{ _{M^{\ot\;n}}}$ for all $n \geq 1$.
\end{lemma}

\bpf
Suppose $m\in M \;,n \in M$ and $a\in \mtr{LKer}_{ _M}\cap \mtr{LKer}_{ _{N}} $. Since $a_1\ot a_2n=a\ot n$ applying $\D\ot \mtr{Id}$ one has $a_1\ot a_2 \ot a_3n=a_1\ot a_2 \ot n$. Thus $a_1\ot a_2m\ot a_3n=a_1\ot a_2m\ot n=a_1\ot m\ot n$.
\epf

\bn{defn}
Let $A$ be a Hopf algebra and $M$ be an $A$-module. We define the Hopf kernel $A_{ _M}$ of $M$ as the largest sub-bialgebra of $A$ contained in ${\bf \mtc{S}}_{ _M}$.
\end{defn}

If $A$ is finite dimensional then $A_{ _M}$ is also the largest Hopf subalgebra of $A$. It is easy to see that in the case of a semisimple Hopf algebra this notion of kernel coincides with the kernel $A_{ _M}$ of the module $M$ introduced in \cite{Bker}.

\bn{lemma}\label{largestsubcoalgera}
Let $A$ be a Hopf algebra and $M$ be an $A$-module. Then:

1) The kernel $A_{ _M}$ is the largest subcoalgebra in ${\bf \mtc{S}}_{ _{M}}$.

2) The kernel $A_{ _M}$ is the largest subcoalgebra in $\mtr{LKer}_{ _{M}}$.

3) The kernel $A_{ _M}$ is the largest subcoalgebra in $\mtr{RKer}_{ _{M}}$.
\end{lemma}

\bpf
1) If $C$ is a subcoalgebra of $A$ contained in ${\bf \mtc{S}}_{ _{M}}$ then $<C>=\oplus_{n \geq 0}C^n$ is a sub-bialgebra of $A$. Thus one has $<C> \subset A_M$.

2)Clearly $A_{ _M}$ is a subcoalgebra of  $\mtr{LKer}_{ _{M}}\subset {\bf \mtc{S}}_{ _{M}}$. On the other hand any subcoalgebra of ${\bf \mtc{S}}_{ _{M}}$ is by definition included in $\mtr{LKer}_{ _{M}}$. Thus the largest subcoalgebra of ${\bf \mtc{S}}_{ _{M}}$ is $A_{ _M}$ and coincides with the largest subcoalgebra of $\mtr{LKer}_{ _{M}}$.

3)The proof of 3) is similar to that of 2).
\epf

\bp\label{Hopfkernel}
Let $A$ be a Hopf algebra and $M$ be an $A$-module. Then
$$A_{ _M}=\{a\in A\;|\; a_1\ot a_2m \ot a_3 =a_1\ot m \ot a_2\;\;\text{for all}\;\;m \in M\}$$
\ep

\bpf
Let $A'_{ _M}=\{a\in A\;|\; a_1\ot a_2m \ot a_3 =a_1\ot m \ot a_2\;\;\text{for all}\;\;m \in M\}$. Clearly $A_{ _M}\subset A'_{ _M}$, since $A_{ _M}$ is a coalgebra. On the other hand it easy to check that $A'_{ _M}$ is a sub-bialgebra of $A$ contained inside $\bf{\mtc{S}}_{ _M}$. Maximality of the kernel $A_{ _M}$ implies the other inclusion.
\epf

Let $\pi:A\ra B$ a Hopf map. Recall that the Hopf kernel of $\pi$ was defined in \cite{AD} as:

\begin{equation}\label{Hopfk}
    \mtr{HKer}(\pi)=\{a\in A\;|\; a_1\ot \pi(a_2) \ot a_3 =a_1\ot \pi(1) \ot a_2\}
\end{equation}

Also the left and right kernels of $\pi$ are defined as $\mtr{LKer}(\pi)=A^{co\;\pi}$ and $\mtr{RKer}(\pi)=^{co\;\pi}A$.

Now, let $M$ be an $A$-module and $I_{ _M}=\cap_{n \geq 0}\mtr{Ann}_A(M^{\ot\;n})$. Then $I_{ _M}$ is a Hopf ideal \cite{PQ} that will be called the Hopf ideal generated by $M$ for the rest of the paper. Also let $\pi_M:A\ra A/I_{ _M}$ be the canonical projection.

\bn{thm}\label{catchr}
Suppose that $M$ is a finite dimensional module over a finite dimensional Hopf algebra $A$. Let $I_{ _M}=\cap_{n \geq 0}\mtr{Ann}_A(M^{\ot\;n})$ and $\pi:A\ra A/I_{ _M}$ be the canonical projection. Then

1) $A^{co\;\pi}=\mtr{LKer}_{ _M}$ and $^{co \;\pi}A= \mtr{RKer}_{ _M}$.

2) $\mtr{Hker}(\pi)=A_M$.
\end{thm}
\bpf

1)If $a \in A^{co\;\pi}$ then $a_1\ot a_2m=a_1\ot\pi(a_2)m=a\ot\pi(1)m=a\ot m$,
therefore $a \in \mtr{LKer}_{ _M}$. Conversely suppose that $a \in \mtr{LKer}_{ _M}$ and let  $\D(a)=\sum_{i=1}^sa_i \ot x_i$ with $a_i$ a $k$-basis of $A$. Then $x_i\in \mtr{LKer}_{M}$ and the Lemma \ref{power} implies that $x_i \in \mtr{LKer}_{ _{M^{\ot\;n}}}$ for all $n \geq 1$. Therefore $x_i-\eps(x_i)1 \in I_{ _M}$ which implies that $\pi(x_i)=\eps(x_i)1$. Thus $A^{co\;\pi}=\mtr{LKer}_{ _M}$ and applying the antipode $S$ it follows that $^{co \;\pi}A= \mtr{RKer}_{ _M}$.

2)It is easy to see that $A_M \subset \mtr{HKer}(\pi)$. On the other hand since $\mtr{HKer}(\pi)$ acts trivially on $M$ it follows that $ \mtr{HKer}(\pi)\subset A_M$ by maximality of $A_{ _M}$.
\epf

The following Corollary can be regarded as a generalization of Brauer's theorem for groups.

\bn{cor}\label{charofim}Suppose that $M$ is a finite dimensional module over a finite dimensional Hopf algebra $A$. Then $$\cap_{n \geq 0}\mtr{Ann}_A(M^{\ot\;n})=\omega(\mtr{LKer}_{ _M})A.$$
\end{cor}

For a coideal subalgebra $S$ of $A$ denote by $\eps_{ _S}$ the character of the left trivial $S$-module. Then $\eps_{ _S}$ is the restriction of the counit $\eps$ to $S$.

\bn{rem}
One has that $A \ot_Sk$, the trivial left $S$-module induced to $A$, is isomorphic to
$A/AS^+$ via the map $a\ot_S1 \mapsto \bar{a}$.
\end{rem}

Next Lemma is the first item of Theorem 1.1 of \cite{Tkq}.
\bn{lemma}\label{coinv}
Suppose that the antipode $S$ of the Hopf algebra $A$ is bijective.
Let $S$ be a coideal subalgebra of $A$ and $\pi:A\ra A/AS^+$ the
canonical coalgebra projection. If $A$ is left $S$-faithfully flat then $A^{co\;\pi}=S$.
\end{lemma}

\bn{rem}\label{2.3}
From Lemma 4.2 of \cite{Tkq} it follows that $S$ is normal whenever $AS^+\subset S^+A$.
\end{rem}

\bn{prop}\label{closure}
Let $L$ be a normal left coideal subalgebra of a Hopf algebra $A$ with bijective antipode $S$. Then

\begin{equation}\label{closu}
    \mtr{LKer}_{ _{\eps_L\uw^{ ^A}_{ _L}}}=L
\end{equation}

\end{prop}

\bpf
We show first that $L \subset \mtr{LKer}_{ _{\eps_{ _L}\uw_L^A}}$. Indeed,
for all $ l \in L$ and $a \in A$ one has
$$l_1\ot l_2(Sa \ot_L1)=l_1\ot Sa_1\ot_L (a_2l_2Sa_3)1)=l \ot (Sa \ot_L 1)$$
since $L$ is normal. Since the antipode $S$ is bijective it follows that $L \subset \mtr{LKer}_{ _{\eps_{ _L}\uw_L^A}}$.

Let $\pi: A \ra A//L$ be the canonical projection. It follows that $A^{co\;\pi}=L$ by Lemma \ref{coinv}. Suppose now that  $a \in \mtr{LKer}_M$ with $M=k_{ _L}\uw^{ ^A}_{ _L}\cong A/AL^+$. Then
$a\ot m=a_1\ot a_2m=a_1\ot \pi(a_2)m$ for all $m \in M$.
In particular for $m=\bar{1}$ one has $a_1\ot \pi(a_2)=a\ot \pi(1)$.
Therefore $a \in A^{co\;\pi}=L$.
\epf

\bn{cor}\label{nonssBrauer}
Suppose that $L$ is a normal coideal subalgebra of $A$ and $M:=k\uw^{ ^A}_{ _L}$ is the trivial $L$-module induced up to $A$. Let $I_{ _M}$ be the Hopf ideal generated by $M$ in $A$ and $\pi:A\ra A/I_{ _M}$ be the canonical projection.
Then $A^{co\pi}=L$, i.e $A/I_{ _M}=A//L$.
\end{cor}

\bpf
By Proposition \ref{closure} it follows that $\mtr{LKer}_{ _M}=L$. On the other hand Corollary \ref{charofim} implies that $A^{co\pi}=\mtr{LKer}_{ _M}$.
\epf

Next we will give two examples of left (right) kernels and kernels.
\bn{example} 1) Let $A$ be a Hopf algebra and consider the left adjoint action of $A$ on itself. Then $L:=\mtr{LKer}(A)$ is the largest central coideal subalgebra of $A$. In this situation the kernel $A_{ _A}$ of the adjoint action coincides with the largest central sub-bialgebra of $A$ and it is in fact a Hopf subalgebra called, the Hopf center of $A$ in \cite{Andr}.

2) Let $A$ be a finite dimensional Hopf algebra and $$L:=\cap_{S \in \mtr{Irr}(A)}\mtr{LKer}(S).$$ Then $A//L$ is the largest semisimple Hopf algebra quotient of $A$. Moreover $A$ has Chevalley property if and only if the Jacobson radical $\mtr{rad}(A)$ equals $\omega(L)A$.
\end{example}

\bn{rem}
\end{rem}
Let $\pi:A \ra H$ be a Hopf algebra map with $A$ and $H$ finite dimensional Hopf algebras. Then $H$ can be regarded as $A$-module via $\pi$ and let $L:=\mtr{LKer}_A(H)$. It is to see that in this case $A//L$ is isomorphic with the Hopf image of $\pi$ as defined in \cite{BB}. Thus $\pi$ is inner faithful if and only if $L$ is trivial.
\subsection{Core of a coideal subalgebra}
Let $S$ be a left coideal subalgebra of $A$. If $L$ and $K$ are left normal coideal subalgebras of $A$ contained in $S$ then it is easy to see that $LK$ is also a left normal coideal subalgebra contained in $S$. Thus one can define $L:=\mtr{core}(S)$ as the largest left normal coideal subalgebra of $A$ contained in $S$.

The next Proposition gives a description of the core of a coideal subalgebra.
It generalizes Theorem 3.7 and Remark 3.8 from \cite{Bker}.

\bp\label{core}
Suppose that $S$ is a coideal subalgebra of $A$ and let
$L:=\mtr{LKer}_{ \eps_{ _S}\uw^{^A}_{ _S}}.$ Then $\mtr{core}(S)=L$.
\ep

\bpf
First we show that $L\subset S$. Since $k\uw^{^A}_{ _S}=A/AS^+$ one has that
\begin{equation}\label{core}
    L=\{l \in A\;|\; l_1\ot \pi(l_2a)=l\ot \pi(a) \;\;\text{for all}\;a\in A\}
\end{equation}

Thus for $a=1$ one gets that $L\subset A^{co\; \pi}=S$. Now suppose that
$K$ is any left normal  coideal subalgebra of $A$ contained in $S$. Then one
has a canonical projection of $A$ modules $A/AK^+\ra A/AS^+$ and Proposition
\ref{seseq} implies $K=\mtr{LKer}_{ _{A/AK^+}}\subset L$.
\epf

\bn{cor}\label{FR}
Let $S$ be a coideal subalgebra of $A$. Then $S$ is normal if and only if $\eps_{ _S}\uw_{ _S}^A\dw^A_{ _S}=\frac{|A|}{|S|}\eps_{ _S}$
\end{cor}

\bpf
If $S$ normal the statement follows from Proposition \ref{closure}. The converse follows from previous Proposition and Proposition \ref{largestcoidsubalg}.
\epf
\subsection{On two endofunctors on $A-\mtr{mod}$ and respectively $S-\mtr{mod}$}

Let $A$ be a Hopf algebra. Then
$A$ is a right comodule subalgebra of $A$ with the usual multiplication and comultiplication given by $\D$.

Let $S$ be a left $A$-comodule subalgebra of $A$, i.e a left coideal subalgebra of $A$. Then for any $A$-module $M$ and any $S$-module $V$ the comodule structure  $\rho: S \ra A \ot S$ defines via pullback, an $S$-module structure on $M \ot V$. Denote this module structure by $M \odot V$.

This makes the category of $S$-modules a left module category over the tensor category $A$-modules.

\bn{prop}\label{mainrelation}
Let $S \subset A$ be a right $A$-comodule subalgebra of $A$. Then $M\ot V\uw^A_S \cong (M\dw^A_{ _S} \odot V)\uw^A_S$ for any $S$-module $V$ and any $A$-module $M$.
\end{prop}

\bn{proof}
The map $T: A\ot_S(M \odot V) \ra M \ot (A\ot_SV)$ given by $a\ot_S(m \ot v) \mapsto a_1 \ot a_1m \ot (a_2 \ot_Sv)$
is a well defined map and a morphism of $A$-modules with inverse given by $m \ot (a\ot_S v ) \mapsto a_2 \ot_S (S^{-1}(a_1)m \ot v)$.
\end{proof}

\bn{cor}\label{tenspwr}
Let $S \subset A$ be a right $A$-comodule subalgebra of $A$. Then

\begin{equation}\label{touse}
M\ot \eps_S\uw^A_S=M\dw_S^A\uw_S^A
\end{equation}

for any $A$-module $S$.
\end{cor}

\bn{proof}
Put $V=k$, the trivial $S$-module
in the above relation. Note also that $M\odot k=M\dw^A_S$.
\end{proof}

Define the endofunctors:

\begin{equation}\label{T}
\mtc{T}: S-\mtr{mod} \ra S-\mtr{mod}\;\;\;\text{given by}\;\;\;\mtc{T}(V)=V\uw_S^A\dw_S^A
\end{equation}

for any $S$-module $V$. Also define

\begin{equation}\label{V}
    \mtc{V}:  A-\mtr{mod} \ra  A-\mtr{mod}\;\;\;\text{given by}\;\;\;\mtc{V}(M)=M\dw_S^A\uw_S^A
\end{equation}

for any $A$-module $M$.

The following Lemma is straightforward. It shows that the restriction functor  $\mtr{res}: A-\mtr{mod}\ra  S-\mtr{mod}$ is a morphism of $A- \mtr{mod}$ categories.

\bl
Let $S$ be a coideal subalgebra of $A$ and $M$, $N$ be two $A$-modules. Then $$(M\ot N)\dw_S^A=M\odot N\dw_S^A$$
\el

Then one has the following relations:

\bn{lemma}
For any $S$-module $V$ and any $A$-module $M$ it follows that:
\bn{enumerate}
\item $$\mtc{V}^n(M)=M\ot(\eps_S\uw_S^A)^n$$
\item $$\mtc{T}^{n+1}(V)=V\uw^A_S\odot \mtc{T}^{n}(\eps_S)$$
\end{enumerate}
for all $n\geq 1$.
\end{lemma}

\bpf
The first statement follows from Corollary \ref{tenspwr}. The second statement is easily proven by induction on $n$.
\epf

In the next Proposition we need the Frobenius-Perron theory on Grothendieck rings of hopf algebras developed in \cite{EO}.

\bp\label{regchar}
Let $A$ be a finite dimensional Hopf algebra and $L:=\mtr{LKer}_{ _M}$ be the kernel of a finite dimensional $A$-module. Then $\eps_{ _L}\uw^A_{ _L}$ is the regular character of $A//L$.
\ep

\bpf
For any left $A//L$-module one has that $N\dw^A_{ _L}=|N|\eps_{ _L}$. Applying previous lemma it follows that $N\otimes \eps_{ _L}\uw^A_{ _L}=|N|\eps_{ _L}\uw^A_{ _L}$. Thus $\eps_{ _L}\uw^A_{ _L}$ is the common Frobenius-Perron eigenvector of the operators of left multiplication by all $A//L$-modules, thus the regular character of $A//L$ (see \cite{EO}).
\epf
\section{Depth two coideal subalgebras}\label{d2}
In this section we work over a commutative ring $R$. Recall that an extension of $R$-algebras
$B \subset A$ is called right (left) depth two if the module $A\ot_BA$ is a direct summand in
$A^n$ in the category $\;_B\mtr{Mod}_A$ (respectively $_A\mtr{Mod}_B$) for an arbitrary $n \geq 1$.

In this section we will show that a left coideal subalgebra of a Hopf algebra $A$ is
right depth two if and only if it is left normal, i.e closed under the left adjoint action.
For this reason in this section we have to work with right modules instead of left. Our treatment is very similar to the one used in \cite{BJ} for depth two Hopf subalgebras.

\subsection{Depth two coideal subalgebras}

\bn{prop}\label{mapbeta2.2, 2.4}
Let $S$ be a right coideal subalgebra of a Hopf algebra $A$.

1)Then the map
\begin{equation*}
\beta:A\ot_SA \rightarrow A\ot A/S^+A\;\;\; \text{given by}\;\;\;a\ot_Sb\mapsto ab_1\ot\bar{b_2}
\end{equation*}

is a well defined morphism of $(A,S)$-bimodules. The $(A,S)$-bimodule actions are
defined as $c(a\ot_S b)x=ca\ot bx$ on $A\ot_SA$ and $c(a\ot \bar{b})x=cax\ot b$ for
$a, b, c \in A$ and $x \in S$.

2) If $AS^+ \subset S^+A$ then the above map $\beta$ is an isomorphism.
\end{prop}

\bpf
The first item follows by direct computation. For the second item consider
\begin{equation*}
\gamma :A\ot A/S^+A\rightarrow A\ot_SA \;\;\; \text{given by}\;\;\;a\ot_S\bar{b}
\mapsto ab_1\ot_SSb_2
\end{equation*}
It is not hard to check that if $AS^+ \subset S^+A$ then $\gamma$ is well defined.
Moreover $\gamma$ is an inverse for $\beta$.
\epf

For the rest of this section let $\bar{A}:=A/AS^+A$ and $\pi:A \ra \bar{A}$ be
the canonical projection. Recall that the extension $A/S$ is $\bar{A}$ -right
Hopf Galois \cite{Montg} if $A^{co\; \pi}=S$ and the canonical map $\beta$ is bijective.

\bp\label{Galois2.5, 2.6}
If $A$ is left or right faithfully flat over $S$ and $AS^+ \subset SA^+$ then
$A/S$ is a right $\bar{A}$-Galois extension.
In particular $S$ is a left normal coideal subalgebra.
\ep

\bpf
By Proposition \ref{mapbeta2.2, 2.4} one knows that $\beta$ is bijective.
If $A$ is left faithfully flat
over $S$ the second theorem from section 13.1 of \cite{wat} implies that $S$ equals the
equalizer of $i_1$ and $i_2$ where $i_1:A\ra A\ot_SA\;\; a\mapsto a\ot_S1$ and
$i_2:A\ra A\ot_SA\;\; a\mapsto 1\ot_Sa$. Since $\beta$ is bijective this coincides with
the equalizer of $i_1\circ \beta$ and $i_2\circ \beta$.
But this last equalizer is exactly $A^{co\;\pi}$. Thus $S$ is also closed under the left
adjoint action.

If $A$ is right faithfully flat
over $S$ a "right version" of the same second theorem from section 13.1 of \cite{wat}
would also imply that $S$ is the equalizer $i_1$ and $i_2$.
\epf

\bt\label{2.8}
If $AS^+ \subset SA^+$ and $A$ is finitely generated projective as left $S$-module then $S$
is of right depth two inside $A$.
\et

\bpf
Using Lemma 2.7 from \cite{BJ} it follows that $\bar{A}|R^n$ in $\mtr{Mod}_R$ and therefore
$A\ot_SA\cong A\ot \bar{A}|A\ot R^n \cong A^n$ in $_A\mtr{Mod}_S$.
\epf

\bn{prop}\label{d2implies2.9} Let $S$ be a left coideal subalgebra of $A$ such that
$A$ is faithfully flat $S$-module. If $S$ has right depth two inside $A$ then $AS^+\subset S^+A$.
\end{prop}

\bpf
If $S$ has right depth two inside $A$ then $k\ot_A(A\ot_SA)$ is a direct summand in
$(k \ot_AA)^n$ in $\mtr{Mod}_S$. Note that $k \ot_SA$ is the trivial $S$-modules and therefore
$k\ot_SA$ divides $k^n$ in $S$-mod. Since $S^+$ annihilates the $S$-module $k$ this implies
that the two sided ideal generated $AS^+A$ annihilates the $A$-module $k \ot_SA$.
Previous remark implies that $AS^+A =S^+A$ and in particular $AS^+ \subset S^+A$.
\epf

\bn{cor}
Let $S$ be a coideal subalgebra of $A$ such that $A$ is faithfully flat over $S$. Then $S$ is a right depth two subalgebra of $A$ if and only if $S$ is normal, i.e closed under the left adjoint action.
\end{cor}
\subsection{Rieffel's normality for coideal subalgebras}
Let $B \subset A$ an extension of finite dimensional $k$-algebras. An ideal $J$ of $B$ is called $A$-invariant if $AJ=JA$. Following \cite{Ri} the extension $A/B$ is called normal
if for every maximal two sided ideal $I$ of $A$  the ideal $B \cap I$ is $A$-invariant.

If $S$ is closed under adjoint action and $A$ has bijective antipode then it easy to verify
that $S$ is normal in Rieffel's sense. Indeed, if $I$ is a two-sided ideal of $A$ and
$x \in S\cap I$ then $ax=a_1xS(a_2)a_3 \in (I\cap S)A$. Also $xSa=a_3Sa_2xSa_1 \in A(I \cap S)$.
Conversely, if a coideal subalgebra $S$ of a semisimple Hopf algebra $A$ is normal
in Rieffel's sense then $AS^+=S^+A$ since $S^+=A^+\cap S$ is a maximal two sided ideal of $S$.
Then Proposition \ref{Galois2.5, 2.6} implies that $S$ is closed under the left adjoint action.
Thus normality, depth two and Rieffel's normality coincide for coideal subalgebras
of a Hopf algebra with bijective antipode.

\bn{rem}\label{reifn}
If $A$ is semisimple it will be shown in the next section that $S$ is also semisimple.
In this case left depth two coincides with right depth two by Theorem 4.6 of \cite{BKK}.
\end{rem}

This remark brings up the question whether left and right depth two coincide on coideal
subalgebras. Further one can ask whether the extension $A/S$ is a quasi-Frobenius extension, at least in the finite dimensional case. Results from \cite{Sk} shows that $S$ is a Frobenius algebra if $S$ is finite dimensional.

\section{The semisimple case}\label{ss}
Throughout of this section we assume that the Hopf algebra $A$ is semisimple. Let $S$ be a left coideal subalgebra of $A$.

Since $A$ is free over $S$ \cite{Sk} there is a decomposition $A=S\oplus R$ as left $S$-modules. Consider $\Lam_{ _A}=x+r$ the decomposition of the idempotent integral $\Lam_{ _A}$ in the
above direct sum. Then clearly $sx=\eps(s)x$ for all $s \in S$ and $\eps(x)=1$.
Similarly since $A$ is free as right $S$ there is a decomposition $A=S\oplus R'$ as
left $S$-module. Consider $\Lam_{ _A}=y+r'$ the decomposition of the idempotent
integral $\Lam_{ _A}$ in the above direct sum. Then clearly $ys=\eps(s)x$ for all
$s \in S$ and $\eps(y)=1$. Thus $x=yx=y$.

\bn{lemma}
Let $S$ be a left coideal subalgebra of a semisimple Hopf algebra $A$. Then $S$ is also semisimple.
\end{lemma}

\bpf
We will use a Maschke type argument for the category of $S$-modules. Suppose that $W$ is an $S$- submodule of $V$ and $f:V\ra W$ a $k$-linear projection. Then it easy to check  $\tilde{f}:V\ra W$  given by
\begin{equation*}
    \tilde{f}(v)=\sum Sx_1f(x_2v)
\end{equation*}

is an $S$-projection. Indeed one has that $\tilde{f}(sv)=\sum Sx_1f(x_2sv)=s_1S(s_2)Sx_1f(x_2s_3v)=s\tilde {f}(v)$ and $\tilde{f}(w)=Sx_1x_2w=w$.
\epf

Thus the element $x$ from above is the central idempotent corresponding to the trivial module $\eps_{ _S}$ of $S$.
\subsection{Rieffel's equivalence relation for a coideal subalgebra}
Let $S \subseteq A$ be a left coideal subalgebra of $A$.

Let $V$ and $W$ two $S$-modules.  that
We say that $V \sim W$ if and only if there is
a simple $A$-module $M$ such that $V $ and $W$ are
both constituents of $M\downarrow_{ _S}^A$. The relation $\sim$
is reflexive and symmetric but not transitive in general.
Its transitive closure is denoted by $\approx$.
Thus we say that $V \approx W$ if and only if there is $m \geq 1$ and a sequence
$V_{i_0}, V_{i_1},\cdots, V_{i_{m-1}}, V_{i_m}$ of simple $S$-modules such that $V=V_{i_0} \sim V_{i_1}\sim V_{i_2} \sim \cdots \sim V_{i_{m-1}} \sim V_{i_m}=W$. As explained in \cite{BKK} it follows that $V\approx W$ if and only if there is $n \geq 0$ such that $V$ is a constituent to $\mtc{T}^n(W)$. This is equivalent to $W$ to be a constituent of $\mtc{T}^n(V)$.

We denote the above equivalence relation by $d_{ _S}^A$. This equivalence relation is considered in \cite{Ri} in the context of any extension of semisimple Hopf algebras.

Similarly one can define an equivalence relation $u_B^A$ on the set of irreducible $A$-modules. We say that $M\sim N$ for two simple $A$-modules $M$ and $N$ if and only if their restriction to $S$ have a common constituent. Then $u_{ _S}^A$ is the transitive closure of $\sim$. Similarly $M$ and $N$ are equivalent if and only if there is $n \geq 0$ such that $M$ is a constituent to $\mtc{V}^n(N)$. This is equivalent to $N$ to be a constituent of $\mtc{V}^n(M)$.

\subsection{Tensor powers of a character}
The following Theorem can be viewed as a generalization of Brauer's theorem for groups:

\bn{thm}\label{Brauer} Let $A$ be a semisimple Hopf algebra and $M$ be an $A$-module
with character $\ch$. If $L:=\mtr{LKer}_{ _{M}}$ then the irreducible modules of $A//L$ are precisely all the  irreducible constituents the tensor powers $M^{\ot\;n}$ with $n \geq 0$.

\end{thm}

\bn{proof} One has that $A//L$ is a semisimple Hopf algebra. On the other hand since $L$
is left normal Corollary \ref{FR} implies that the simple $A$-submodules of
$k\uw_L^A$ and the simple $A$ -submodules of $A//L$ coincide. Then description of the Hopf ideal $I_{ _M}$ given in Corollary \ref{charofim} implies the conclusion.
\end{proof}

Next Corollary follows from Proposition 3.3 from \cite{Bker}.

\bn{cor}\label{central} Let $A$ be a semisimple Hopf algebra and $M$ an $A$ module with character $\ch$. If $\ch$ is central in $A^*$ then $\mtr{Lker}_{ _{\ch}}=\mtr{Rker}_{ _{\ch}}=A_{ _{\ch}}$.
\end{cor}

\bpf
Let $L:=\mtr{Lker}_{ _{\ch}}$. Then from the previous Theorem and Proposition 3.3 of \cite{Bker} the quotient Hopf algebras $A//L$ and $A//A_{ _{\ch}}$ have the same irreducible representations, namely the irreducible representations of all tensor powers of $M$. It follows that $\dim_kL=\dim_kA_{ _{\ch}}$. Since $A_{ _{\ch}}\subset L$ then one has $A_{ _{\ch}}=L$. Then $A_{ _{\ch}}=\mtr{Rker}_{ _{\ch}}$ by applying the antipode $S$.
\epf

In view of the last corollary the next theorem can be viewed as an extension of Corollary 6.5 for group algebras from \cite{BKK}.

\bn{thm}\label{coincidenceofequivclss}Let $A$ be a semisimple Hopf algebra and $S$ be a
coideal subalgebra of $A$. Let $L:=\mtr{core}(S)$. Then the equivalence relations
$u_{ _S}^A$ and $u_{ _L}^A$ coincide.
\end{thm}

\bpf
From Proposition \ref{core} it follows that $L:=\mtr{LKer}_{\eps_S\uw^A_S}$.
Two $A$-modules $M$ and $N$ are equivalent if and only if $N$ is a constituent of
$\mtc{V}^n(M)$ for some $n \geq 1$. From the formula of $\mtc{V}^n(M)$ it follows that
$M$ and $N$ are equivalent if and only if $N$ is a constituent of $M\ot (\eps_S\uw^A_S)^n$
for some $n \geq 1$.

But from Theorem \ref{Brauer} it follows that $\eps_{L}\uw^A_L$ has as
constituents all the irreducible constituents of all tensor powers  $(\eps_S\uw^A_S)^n$
with $n \geq 0$. Thus two $A$-modules $M$ and $N$ are equivalent relative to $S$ if and only if
they are equivalent relative to $L$.
\epf
\subsection{The integral element of $S$}
Let $x_{ _S}:=x$ denote the element from the beginning of this section.

In this situation one has $AS^+=A(1-x_{ _S})$. Indeed if $s \in S^+$ then $s=sx_{ _S}+s(1-x_{ _S})=\eps(s)x_{ _S}+s(1-x_{ _S})=s(1-x_{ _S})$.

\bp\label{car} Let $A$ be a semisimple Hopf algebra and $S$ be a coideal subalgebra of $A$. Then the following statements are equivalent:

1)$S$ is normal in $A$

2)$\eps_{ _S}$ by itself from an equivalence class of $d^A_{ _S}$.

3)The element $x_{ _S}$ is central in $A$.
\ep

\bpf
Corollary \ref{FR} implies the equivalence between the first two items. Proposition 3.1 from \cite{BKK} shows the equivalence between the last two items.
By results from \cite{Tkq} it follows that $S$ is normal if and only if $AS^+=S^+A$.
\epf

Note that the equivalence between first and third item generalizes Proposition 1 from \cite{Masnr}

\bibliographystyle{plain}
\bibliography{cidalgs}
\end{document}